\documentclass{article}
\usepackage{graphicx} 
\usepackage{amsmath}

\newcounter{local}
\newcounter{locallocal}
\newcommand{\scl}{\stepcounter{local}}
\newcommand{\h}{\hspace{1cm}}

\newcommand{\by}{\begin{eqnarray}}
\newcommand{\ey}{\end{eqnarray}}
\newcommand{\bys}{\begin{eqnarray*}}
\newcommand{\eys}{\end{eqnarray*}}

\begin{document}

\begin{center}
{\bf On a Modified Mathematical Model Arising from a Trojan Y Chromosome Strategy}
\end{center}

\vspace{1cm}
\begin{center}
Hong-Ming Yin\footnote{Email address: hyin@wsu.edu}\\
Department of Mathematics and Statistics \\
Washington State University\\
Pullman, WA 99164, USA.
\end{center}

\ \\
{\bf Abstract}
In this paper, we investigate a modified Trojan Y chromosome (TYC) strategy aimed at eradicating invasive species from natural habitats. The proposed mathematical model enhances the original TYC framework by ensuring the non-negativity of population densities and preventing potential solution blow-up. The new model is formulated as a strongly coupled reaction-diffusion system with distinct diffusion coefficients for each species. We first establish the global well-posedness of the system. Subsequently, a stability analysis is conducted. In particular, we demonstrate that the population densities converge to zero when the birth rate for each species falls below a critical threshold. Additionally, we prove the existence of a positive steady-state solution even as the artificially introduced YY-female population density tends to zero as $t$ tends to $\infty$. Furthermore, we identify a bifurcation in the system's solutions as the birth rate crosses the critical value.

\vspace{1cm}

\vspace{1cm}

\newpage
\section{Introduction}

 Due to globalization, many exotic species are invading a natural habitat. The new species may destroy the natural balance in the original environment and threaten the life of the original species; see \cite{ALNM2006,SK1997,SL2015} for examples. To eliminate this exotic species, some biological researchers introduced a new method by reversing the gene of a male of the invasive species to create an artificial female; see \cite{FFY2004,CPB1999}. These sex-reversing females become YY-supermales who attract females of the invasive species. Due to the reverse of the gene, YY-supermale  and YY female will not be able to produce their next generation, which in turn will eliminate the invading species. This novel approach provides a promising way to reduce the population of invasive species and protect the original ecosystem. 
 
 To describe this strategy, some researchers proposed a mathematical model governed by a system of ODEs (ODE model); see, for examples, \cite{ALNM2006,BPH2020,FFY2004,Wang2014} or by a reaction-diffusion system (PDE model), see, for examples, \cite{GT2006,PG2011,P2011}.  For the reader's convenience, we recall the PDE model here. The population density for each species is given in the following table.

\ \\
\begin{tabular}{|l|l|}
\hline
 f &  XX  \mbox{ female population density (natural femal)}  \\
 \hline
 m &  XY \mbox{male population density (natural male)} \\
 \hline
 s &  YY  \mbox{ supermale population density (produced through mating $r$)}  \\
 \hline
 r &  YY \mbox{  sex-reversed female density (artificially introduced) }\\
 \hline
\end{tabular}

 \ \\
The PDE model is governed by the following reaction-diffusion system (\cite{GT2006}):
\begin{align*}
\frac{\partial f}{\partial t} =& D\Delta f+\frac{\beta}{2}fmg-df,\\
\frac{\partial m}{\partial t} =& D\Delta m+\beta\left(
\frac{1}{2}fm+\frac{1}{2}rm+fs\right)g-d m,\\
\frac{\partial s}{\partial t} =& D\Delta s+\beta\left(\frac{1}{2}rm+rs\right)g-d s,\\
\frac{\partial r}{\partial t} =& D\Delta r+\mu -dr,
\end{align*}
 where
\[ g:=g(f,m,s,r)=1-\frac{f+m+r+s}{K},\]
$\beta$ is the birth rate for each species, $D$ represents the diffusion coefficient and $d$ represents the death rate, the constant $\mu$ is the rate introduced artificially by sex-reversed female, the maximum carrying capacity of each species is denoted by $K$.
 
 Many researchers have studied the above reaction-diffusion system; see \cite{GT2006,HMP2004,GHPT2012,P2011,PG2010} and the references therein. Under certain conditions, a global well-posedness to the reaction-diffusion system is established. In addition, some asymptotic analysis is also performed (\cite{PG2011}). In a recent paper \cite{TBGBP2022}, Takyi et al. constructed some initial values for the Trojan Y chromosome model such that the population density for the YY male becomes negative and the $L^p$-norm for the solution of the system will blow up in finite time if the initial data are sufficiently large. However, these solutions do not make any sense in biological science because the density cannot become negative. This motivates us to make certain modifications to the original PDE model to prevent these unrealistic solutions from occurring. 
 
 In this paper, we modified the above PDE model to prevent some unrealistic solutions that may occur to the original PDE system.  In our modified PDE model, we assume that the diffusion coefficient and the death rate for each species are different. More importantly, we replace the rate $\mu$ of YY-male artificially introduced by $\mu r g$. This modification will ensure that the density of the sex-inversed female population does not exceed the maximum carrying capacity $K$. Let $\Omega\subset R^n$ be a bounded domain with $C^1$-boundary and $Q_T=\Omega\times (0, T]$ for any $T>0.$ We denote by $Q=\Omega\times (0, \infty)$ when $T=\infty$. By assuming a logistic growth for each species, we see that $f, m, s$ and $r$ satisfies the following reaction-diffusion system in $Q$:

\setcounter{section}{1}
\setcounter{local}{1}
 \begin{eqnarray}
  & \frac{\partial f}{\partial t} - \nabla\cdot[a_1(x,t)\nabla f]=& \frac{\beta}{2}fmg-d_1f,\\
 & \frac{\partial m}{\partial t} - \nabla\cdot[a_2(x,t)\nabla f]=& \beta\left(
 \frac{1}{2}fm+\frac{1}{2}rm+fs\right)g^+-d_2m, \scl \\
 & \frac{\partial s}{\partial t} - \nabla\cdot[a_3(x,t)\nabla s]=& \beta\left(\frac{1}{2}rm+rs\right)g^+-d_3 s,\scl \\
 & \frac{\partial r}{\partial t} -\nabla\cdot[a_4(x,t)\nabla r]=&\mu rg-d_4r, \scl
 \end{eqnarray}
where
\[ g:=g(f,m,s,r)=1-\frac{f+m+r+s}{K},\,  g^+=\max\{g, 0\}, \]
$a_i$ and $d_i$ represent the diffusion coefficient and natural death rate for each species, respectively. 

To complete the mathematical model,  the following initial and boundary conditions are imposed:
\by
 & & (f(x,0),,m(x,0),s(x,0),r(x,0))=(f_0(x),m_0(x),s_0(x),r_0(x)),  x\in \Omega,\nonumber\\
 \scl\\
 & & (a_1\nabla_{\nu}f, a_2\nabla_{\nu}m,a_3\nabla_{\nu}s,a_4\nabla_{\nu}r)=0, \, (x,t)\in\partial \Omega\times (0, \infty).\scl
\ey
where $\nabla_{\nu}$ represents the normal derivative on $\partial \Omega.$

 In the above modified model, we replace the original $g$ by $g^+$ in Eq.(1.2) and Eq.(1.3).
 From a biological point of view, $g^+$ simply means that no additional $XY$-male and $YY$-supermale are allowed in the ecological environment when the total population reaches the maximum carrying capacity. We would like to point out that even with this modification the "Quasi-positivity condition" (see \cite{P2010,PS2000}) does not hold for Eq. (1.3) and Eq.(1.4). However, we will prove that the population density for each species for the system (1.1)-(1.4) will be non-negative as long as initial data are non-negative. We would like to make another point that the global existence results obtained for the reaction-diffusion system in \cite{FMTY2021,MT2020} cannot be applied here, since the intermediate entropy condition is not satisfied for system (1.1)-(1.4). In this paper, we use the energy method to derive an a priori bound globally for all species. With this bound, we are able to establish the global existence and uniqueness. Moreover, using ideas from \cite{Temam,YCW2017,Yin2024}, we performed the stability analysis for the constant steady-state solution to the system. We will also prove that there exists a bifurcation when the birth rate $\beta$ changes across the critical value $\beta_0$ (see the definition of $\beta_0$ in Section 2). To the best of our knowledge, these results are new.

The paper is organized as follows. In Section 2 we state the basic assumptions and state the main results. In Section 3, we derive some $a\, priori$ estimates and prove the global existence by employing a Leray-Schauder's fixed-point theorem. In Section 4 we perform the stability analysis for the steady-state solution of system (1.1)-(1.6). Some concluding remarks are given in Section 5.

\section{Preliminary and the Statement of Main Results}

Throughout this paper, the following basic conditions are assumed.\\
H(2.1). Let $a_i(x,t)\in L^{\infty}(Q)$. There exists a positive constant $a_0$ such that
\[ 0< a_0\leq a_i(x,t)\leq \frac{1}{a_0}, \h \forall (x,t) \in Q.\]
H(2.2). Let $d_i(\xi_1,\xi_2,\xi_3,\xi_4)$ be locally Lipschitz continuous and $\mu(x,t)\in L^{\infty}(Q)$. There exist positive constants $D_0>0$ and $D_1$ such that
\[0<D_0 \leq d_i(\xi_1,\xi_2,\xi_3,\xi_4), \mu(x,t) \leq D_1, \forall (\xi_1, \xi_2, \xi_3, \xi_4)\in (R^+)^4, (x,t)\in Q.\]
H(2.3). The initial values are bounded and nonnegative in $\Omega$: 
\[ 0\leq  f_0(x), m_0(x), s_0(x), r_0(x)\leq K, \forall x\in \Omega. \]
Moreover,
\[ ||f_0||_{L^{\infty}(\Omega)}+||m_0||_{L^{\infty}(\Omega)}+||s_0||_{L^{\infty}(\Omega)}
+||r_0||_{L^{\infty}(\Omega)}\leq K.\]

For convenience, we denote by $F_1, F_2, F_3, F_4$ the right-hand side of Eq. (1.1), Eq.(1.2), Eq.(1.3) Eq.(1.4), respectively. The inner product for $L^2(\Omega)$ is denoted by
\[ (u,v)=\int_{\Omega} uvdx, \forall u, v\in L^2(\Omega).\] 
The Sobolev space $H^1(\Omega)$ is standard and the Banach space $V_2(Q_{T}):=L^2(0, T; H^1(\Omega))\bigcap C([0,T];L^2(\Omega))$ (\cite{EVANS}). For a Banach space $X$, its dual space is denoted by $X^*$ and the pair $<u,v>:=u(v)\in R^1, v\in X, u\in X^*.$

The definition of a weak solution to system (1.1)-(1.6) is standard as in \cite{EVANS}. The dual space of $H^1(\Omega)$ is denoted by $H^*(\Omega)$.

\ \\
{\bf Definition.} We say a vector function 
\[ {\bf Z}(x,t):=(Z_1, Z_2, Z_3, Z_4):=(f(x,t),m(x,t),s(x,t),r(x,t))\] is a weak solution to the system (1.1)-(1.6) if 
\[ Z_i(x,t) \in V_2(Q_{T})\bigcap L^{\infty}(Q_{T}), Z_{it}\in H^*(\Omega),\] satisfies
\begin{eqnarray*}
 & &  <Z_{it}.\phi>+(a_i\nabla Z_i, \nabla \phi)=(F_i, \phi), a.e. t\in (0,T]; i=1,2,3,4
\end{eqnarray*}
for any  $ \phi\in H^1(\Omega))$.
Moreover,
\[ \lim_{t\rightarrow 0^+}||{\bf Z}(x,t)-{\bf Z}_0(x)||_{L^{2}(\Omega)^4}=0.\]

\ \\
{\bf Theorem 2.1.} Under the assumptions H(2.1)-(2.3), system (1.1)-(1.6) has a unique weak solution $(f,m,s,r)$, and the weak solution is nonnegative in $Q$. Moreover,
the solution is H\"older continuous in $\bar{Q}\bigcap \{(x,t): t>0\}.$
Furthermore, the weak solution continuously depends on initial data: For any $T>0$,
\[ \sup_{0\leq t\leq T}||{\bf Z}^*-{\bf Z}^{**}||\leq C(T) ||{\bf Z}_{0}^{*}-{\bf Z}_{0}^{**}||_{L^2(\Omega)},\]
where ${\bf Z}^*(x,t)$ and ${\bf Z}^{**}(x,t)$ are the weak solution of (1.1)-(1.6) corresponding to initial data ${\bf Z}_{0}^{*}(x)$ and ${\bf Z}_{0}^{**}(x)$, respectively.

To study the asymptotic behavior of the solution, we assume an additional condition for $\mu(x,t)$.

\ \\
H(2.4). Assume that
\[ \lim_{t\rightarrow \infty}||\mu||_{L^2(\Omega)}=0.\]

Define
\[ \beta_0:=\frac{8(d_1+d_2)}{K}.\]
\ \\
{\bf Theorem 2.2.}  Under the assumptions H(2.1)-(2.4), the weak solution $(f,m,s,r)$ of the system (1.1)-(1.6) has the following asymptotic properties as $t\rightarrow \infty$:
\ \\
(a) If 
\[ 0<\beta<\beta_0,\]
then, there exists only one steady-state solution $(f,m,s,r)=(0,0,0,0)$ to the system (1.1)-(1.6).
Moreover, the steady-state solution $(0,0,0,0)$ is asymptotically stable.

\ \\
(b) If 
\[ \beta=\beta_0,\]
then, there exist two steady-state solutions to the system (1.1)-(1.6):
\[ (f_1,m_1,s_1,r_1)=(0,0,0,0), \, (f_2,m_2,s_2,r_2)=(\frac{4d_2}{\beta}, \frac{4d_1}{\beta},0,0).\]
Moreover, the steady-state solution $(f_1,m_1,s_1,r_1)$ is asymptotically stable, but
$(f_2,m_2,s_2,r_2)$ is unstable.

\ \\
(c) If 
\[ \beta>\beta_0,\]
then, there exist three steady-state solutions $(f_i,m_i,s_i,r_i), i=1,2,3,$ to the steady-state system (1.1)-(1.6), where $s_i=r_i=0, i=1,2,3.$\\
Moreover, $(f_1,m_1, s_1, r_i)$ is always asymptotically stable, $(f_3, m_3, s_3, r_3)$ is asymptotically stable if
 \[ \beta <\frac{8(d_1+d_2)}{K(1-b)^2}\]
 and unstable if 
\[ \beta > \frac{8(d_1+d_2)}{K(1-b)^2},\]
 where
 \[ b:=\sqrt{1-\frac{8(d_1+d_2)}{K\beta}}\]
but the steady-state solution $(f_2,m_2,s_2,r_2)$ is always asymptotically unstable.

Furthermore, a bifurcation occurs for the solution of the system (1.1)-(1.6) when $t$ tends to $\infty$ and $\beta$ changes across the critical value $\beta_0$.

\section{Proof of the Theorem 2.1}

We begin with a comparison lemma for a coupled reaction-diffusion system.

\ \\
{\bf Lemma 3.1.} Let $u$ and $v$ be a weak solutions of the following parabolic system:
\bys
   & & u_t=\nabla \cdot[a_1(x,t)\nabla u]+b_{11}(x,t)u+b_{12}(x,t)v,\h (x,t)\in Q,\\
   & &  v_t=\nabla \cdot [a_2(x,t)\nabla v]+b_{21}(x,t)u+b_{22}(x,t)v, \h (x,t)\in Q,\\
   & &  (a_1\nabla_{\nu}u, a_2\nabla_{\nu}v)=0, \h (x,t)\in \partial \Omega\times (0, \infty),\\
   & & (u(x,0), v(x,0))=(u_0(x), v_0(x)), \h x\in \Omega.
\eys
Assume that 
\[ 0<a_0\leq a_1(x,t), a_2(x,t)\leq A_0<\infty, \h \forall (x,t)\in Q. \]
Suppose $b_{ij}\in L^{\infty}(Q)$ and $ b_{12}(x,t), b_{21}(x,t)\geq 0$ in $Q$.
If $u_0(x)\geq 0, v_0(x)\geq 0$ in $\Omega$, 
then
\[ u(x,t), v(x,t)\geq 0, \h (x,t)\in Q.\]
{\bf Proof.} Let $T>0$ be arbitrary. 
We first assume $u_0(x), v_0(x) \geq \varepsilon >0$ over $\Omega$ for a small constant $\varepsilon>0$.
From the continuity of the solution, we see that $u(x,t), v(x,t)>0$ in $\bar{Q}_T$ for some small $T>0$.
Define
\[ T^*=\sup\{T: u(x,t)>0, v(x,t)>0, \forall (x,t)\in \bar{Q}_{T}\}.\]
If $T^*=\infty$, the proof is done. If $T^*<\infty,$ then at least one of $u(x,t)$ and $v(x,t)$
attains its minimum $0$ at some point $(x^*, T^*)\in \bar{Q}_{T}$. Suppose $u(x^*,T^*)=0$. Hopf's lemma implies that
$(x^*,T^*)$ must be an interior point of $Q_{T^*}$.
By the definition of $T^*$, we see $v(x,t)\geq 0$ in $Q_{T^*}$. Hence, 
\[ d_{12}(x^*,T^*)v(x^*,T^*)\geq 0.\]
This is a contradiction to the strong maximum principle.
Therefore,  $T^*=\infty$. Now we just take an approximation for initial data and take a limit as $\varepsilon \rightarrow 0$ to conclude the desired result.

\hfill Q.E.D.

\ \\
{\bf Lemma 3.2}. (Nonnegativity and boundedness) Under the assumptions H(2.1)-H(2.3) the concentration of each species is nonnegative and bounded:
\[ 0\leq f(x,t), m(x,t), s(x,t), r(x,t)\leq K, \h \forall (x,t)\in Q.\]
{\bf Proof.} Since $\mu(x,t)\geq 0$ and initial values are nonnegative, the maximum principle implies that
\[ f(x,t) \geq 0, r(x,t)\geq 0, \h \forall (x,t)\in Q.\]
To show the nonnegativity of $m$ and $s$, we see that $m(x,t)$
and $s(x,t)$ satisfy the parabolic system in Lemma 3.1 with
\begin{align*}
& b_{11}=\beta\left(\frac{1}{2}f+\frac{1}{2}r\right)g^+-d_2, & b_{12}=\beta fg^+ \geq 0,\\
& b_{21}=\beta g^+ \left(\frac{1}{2}r\right)\geq 0, & b_{22}=\beta g^+r-d_3
\end{align*}
It follows that 
\[ m(x,t) \geq 0, s(x,t)\geq 0, \h \forall (x,t)\in Q.\]

Next, we derive an upper bound for $f, m, s, r$. 

We use the energy method to see
\begin{align*}
& \frac{1}{2}\frac{d}{dt}\int_{\Omega}((f-K)^{+})^2dx +\int_{\Omega}a_1|\nabla (f-K)^{+}|^2dx
+\int_{\Omega}d_1f(f-K)^+dx\\
& =\frac{\beta}{2}\int_{\Omega}fmg (f-K)^{+}dx \leq 0.
\end{align*}
It follows for any $T>0$ that 
\[ (f-K)^+=0, \h (x,t)\in Q_{T}. \] i.e.
\[ 0\leq f(x,t)\leq K, \h (x,t)\in Q.\]
Similarly, from Eq.(1.2)-(1.4) we see
\[ \int_{\Omega}g^+(g-K)^+dx = 0; \h \int_{\Omega}\mu rg(r-K)^+dx \leq 0.\]
It follows that
\[ 0\leq m(x,t), s(x,t), r(x,t)\leq K, \h (x,t)\in Q.\]

\hfill Q.E.D.

For any $\delta>0$, define $Q_{T}(\delta):=\Omega\times (\delta,T].$
\ \\
{\bf Lemma 3.3.} There exist a $\alpha\in (0, 1)$ and a constant $C$ such that
\[ ||f||_{C^{\alpha, \frac{\alpha}{2}}(\bar{Q}_{T}(\delta))}+||m||_{C^{\alpha, \frac{\alpha}{2}}(\bar{Q}_{T}(\delta))}+||s||_{C^{\alpha, \frac{\alpha}{2}}(\bar{Q}_{T}(\delta))}+||r||_{C^{\alpha, \frac{\alpha}{2}}(\bar{Q}_{T}(\delta))}\leq C(\delta),\]
where $C(\delta)$ depends only on known data, $\delta$ and $T$.\\
{\bf Proof.} Since all terms on the right-hand side of (1.1)-(1.4) are bounded from Lemma 3.2, we can apply DiGiorgi-Nash's estimate (\cite{LSU}) for the weak solution of (1.1)-(1.6) to obtain the desired estimate.

\hfill Q.E.D. 

\ \\
{\bf Proof of Theorem 2.1.}   With the $a priori$ estimates in Lemma 3.2 and Lemma 3.3, we can use the standard Leray-Schauder's fixed point theorem to establish the existence. We give only an outline of the proof here.\\
Without loss of generality, we may assume that
all initial data belong to $C^{\alpha, \frac{\alpha}{2}}(\bar{\Omega})$. Otherwise, we can use the local existence and DiGiorgi-Nash theory to obtain a weak solution in $Q_{t_0}$ for a sufficiently small $t_0>0$. Then we use $(f(x,t_0), m(x,t_0), s(x,t_0), r(x,t_0))$ as an initial value to study problem (1.1)-(1.6) in $\Omega\times [t_0, T]$. Let $T>0$ be any fixed number.

 Choose a convex connect set
\[ X=\{(f,m,s,r)\in L^{\infty}(Q_{T})^4: 0\leq f, m, s, r\leq K, \forall (x,t)\in Q_{T}.\}\subset L^{\infty}(Q_{T})^4.\]
\ \\
Step 1: Let $\sigma\in (0,1]$ be a parameter. Define a mapping $M_{\sigma}$ as follows:

Given any ${\bf Z}:=(f,m,s,r)\in X$, define a mapping $M$ from $X$ to $L^{\infty}(Q_{T})^4$ as follows:
\[ M_{\sigma}[{\bf Z}]:=(f^*,m^*,s^*,r^*), \]
where $(f^*, m^*, s^*,r^*)$ is a solution of the following linear system:
\[ (Z_i)_t-\nabla\cdot(a_i\nabla Z_i)=\sigma F_i(f,m,s,r), \h (x,t)\in Q_{T},\]
subjection to initial and boundary conditions:
\[ Z_i(x,0)=\sigma Z_{i0}(x), \, x\in \Omega,\,  a_i\nabla_{\nu}Z_i(x,t)=0, \forall (x,t)\in \partial \Omega\times (0,T],i=1,2,3,4.\]
The standard theory for parabolic system (\cite{EVANS}) implies that the linear system has a unique weak solution $(f^*,m^*,s^*,r^*)$ for any given $(f,m,s,r)\in X$. Hence, the mapping $M_{\sigma}$ is well defined. These a priori estimates in Lemma 3.2 imply that the mapping $M_{\sigma}$ is from $X$ to $X$ for all fixed points of $M_{\sigma}$:
\[ M_{\sigma}[{\bf Z}]:=(f,m,s,r), \]
Moreover, Lemma 3.3 implies that
$M_{\sigma}[{\bf Z}]\subset C^{\alpha, \frac{\alpha}{2}}(\bar{Q}_{T})$.
On the other hand, the embedding operator from $C^{\alpha, \frac{\alpha}{2}}(\bar{Q}_{T})$ to $L^{\infty}(Q_{T})$ is compact. It follows that $M_{\sigma}$ is a compact mapping from $X$ to $X$ and $M_0({\bf Z})=0$. The proof for the continuity of $M_{\sigma}$ is similar to the proof for the continuous dependence in the following. By Leray-Schauder's fixed-point theorem, the mapping $M_{\sigma}$ has a fixed point. This fixed point for $\sigma=1$ is a solution to the original problem (1.1)-(1.6).

\ \\
Step 2. Prove for the continuous dependence and uniqueness in $Q_{T}$.

 Suppose that ${\bf Z}^*(x,t)$ and ${\bf {Z}}^{**}(x,t)$ are two weak solutions to the system (1.1)-(1.6)
corresponding to the initial data ${\bf Z}_0^*(x)$ and ${\bf {Z}}_0^{**}(x)$, respectively. Let
\[ {\bf Z}={\bf Z}^*(x,t)-{\bf {Z}^{**} }(x,t), \h (x,t)\in Q.\]
Let
\[ g^*=1-\frac{f^*+m^*+s^*+r^*}{K}, \h g^{**}=1-\frac{f^{**}+m^{**}+s^{**}+r^{**}}{K}.\]
Obviously,
\[|g^*-g^{**}|\leq \frac{|{\bf Z}|}{K}.\]
Since $d_i(f,m,s,r)$ is locally Lipschitz continuous for all $i=1,2,3,4$, we use the energy method to obtain
\bys
& & \frac{d}{dt}\int_{\Omega}|{\bf Z}|^2dx \leq C\int_{\Omega}|(g^*)^+-g^{**}|(|Z_2|+|Z_3|)dx+C\int_{\Omega}|{\bf Z}|^2dx.
\eys
Note that
\[ y^+=\frac{|y|+y}{2}, \h \forall y\in R^1.\]
It follows that
\bys
& & \big{|}(g^*)^+-g^{**}\big{|} =\big{|}\frac{(|g^*|+g^*)-(|g^{**}|+g^{**})}{2}\big{|}\\
& & \leq \frac{|g^*-g^{**}|+|(|g^*|-|g^{**}|)|}{2}\\
& & \leq |g^*-g^{**}|\\
& & \leq \frac{|{\bf Z}|}{K},
\eys
where at the third step we have used an elementary inequality:
\[\big{|}|a|-|b|\big{|}\leq |a-b|, \h \forall a, b\in R^1.\]
It follows that
\bys
& & \frac{d}{dt}\int_{\Omega}|{\bf Z}|^2dx \\
& & \leq  C\int_{\Omega}|{\bf Z}|(|Z_2|+|Z_3|)dx+C\int_{\Omega}\int_{\Omega} |{\bf Z}|^2dx\\
& & \leq C\int_{\Omega}|{\bf Z}|^2dx,
\eys
where $C$ depends only on known data.

Gronwall's inequality yields
\[ \sup_{0\leq t\leq T}||{\bf Z}||_{L^{2}(\Omega)}\leq C(T)||{\bf Z}_0^*-{\bf Z}_0^{**}||_{L^2(\Omega)},\]
where $C(T)$ depends only on known data and $T$.

The uniqueness follows immediately. This concludes the proof of Theorem 2.1.

\hfill Q.E.D.

\section{ Asymptotic Analysis}

In this section, we analyze the stability for the constant solution of the system (1.1)-(1.4).

We define a $\Omega$-limit set (see\cite{Temam}):
\[W:=\{ \bar{{\bf Z}}:=(\bar{f},\bar{m},\bar{s},\bar{r}): \exists t_n, \lim_{t_n\rightarrow \infty}||{\bf Z}-{\bf \bar{Z}}||_{L^2(\Omega)}=0\}.\]

\ \\
{\bf Lemma 4.1.} Suppose that Assumption H(2.4) is valid.
Then, 
\[ \lim_{t\rightarrow \infty} ||r||_{L^2(\Omega)}=\lim_{t\rightarrow \infty}||s||_{L^2(\Omega)}=0.\]
{\bf Proof.} Since the solution $(f,m,s,r)$ is uniformly bounded, we see
\[ \frac{1}{2} \frac{d}{dt}\int_{\Omega}r^2dx+a_0\int_{\Omega}|\nabla r|^2dx +d_4\int_{\Omega}r^2dx
\leq C\int_{\Omega}\mu^2dx.\]
Gronwall's equation implies
\[ \lim_{t\rightarrow \infty}||r||_{L^2(\Omega)}=0.\]
Similarly, from Eq.(1.3) we see
\[  \frac{1}{2} \frac{d}{dt}\int_{\Omega}s^2dx+a_0\int_{\Omega}|\nabla s|^2dx +d_3\int_{\Omega}s^2dx
\leq C\int_{\Omega}r^2dx.\]
We use Gronwall's equality again to obtain
\[ \lim_{t\rightarrow \infty} \int_{\Omega}s^2dx=0.\]

\hfill Q.E.D.

\ \\
{\bf Proof of Theorem 2.2.}

The results in Lemma 4.1 implies that both $r(x,t)$ and $s(x,t)$ converge to $0$ uniformly in $L^2(\Omega)$ as 
$t\rightarrow \infty$, we can neglect the term $ \frac{1}{2}rm + fs$ in Eq.(1.2) in the investigation of the dynamics of the solution ${\bf Z}(x,t)$. Moreover, we will see that all steady-state solutions satisfy
\[ f+m+s+r\leq K.\]
Hence, 
\[ g^+(f,m,s,r)=g(f,m,s,r).\]

Consider the following system
\setcounter{section}{4}
\setcounter{local}{1}
 \begin{eqnarray}
 & & \frac{\partial f}{\partial t} = \nabla\cdot[a_1\nabla f] f+\frac{\beta}{2}fmg_1-d_1f,\\
 & & \frac{\partial m}{\partial t} = \nabla \cdot[a_2\nabla m]+
 \frac{\beta}{2}fmg_1-d_2m, \scl
\end{eqnarray}
where
\[ g_1=1-\frac{f+m}{K}.\]

Now we can easily find the constant solution for the following system:
\bys
& & \frac{\beta}{2}fmg_1-d_1f=0,\\
& &  \frac{\beta}{2}fmg_1-d_2m=0.
\eys

\ \\
Case 1. If 
\[ \beta<\beta_0:=\frac{8(d_1+d_2)}{K}, \]
then $(f^*, m^*)=(0,0)$ is only constant solution for the steady-state system (4.1)-(4.2).
No bifurcation occurs. The  $\omega$-limit set contains only one point:
\[ G=\{(0,0)\},\]
which is a global attractor (\cite{Temam}).

\ \\
Case 2. If 
\[ \beta=\frac{8(d_1+d_2)}{K},\]
then there exist two steady-state solutions:
\[ (f_1,m_1)=(0,0), \h (f_2, m_2)=(\frac{4d_2}{\beta}, \frac{4d_1}{\beta}). \]

Case 3. If 
\[ \beta>\frac{8(d_1+d_2)}{K},\]
then, there exist three steady-state solutions:
\[ (f_1,m_1)=(0,0), \, (f_2, m_2)=(\frac{Kd_2(1+b)}{2(d_1+d_2)}, \frac{Kd_1(1+b)}{2(d_1+d_2)}), \, (f_3, m_3)
=(\frac{Kd_2(1-b)}{2(d_1+d_2)}, \frac{Kd_1(1-b)}{2(d_1+d_2)}), \]
where
\[ b:=\sqrt{1-\frac{8(d_1+d_2)}{K\beta}}\]
Define 
\begin{eqnarray*}
& & F_1:=\frac{\beta}{2}fmg_1-d_1f\\
& & F_2:=\frac{\beta}{2}fmg_1-d_2m
\end{eqnarray*}
A direct calculation yields
\begin{eqnarray*}
& & \frac{\partial F_1}{\partial f}=\frac{\beta m}{2}(g_1-\frac{f}{K})-d_1,\\
& & \frac{\partial F_1}{\partial m}=\frac{\beta f}{2}(g_1-\frac{m}{K}),\\
& & \frac{\partial F_2}{\partial f}=\frac{\beta m}{2}(g_1-\frac{f}{K}),\\
& & \frac{\partial F_2}{\partial m}=\frac{\beta f}{2}(g_1-\frac{m}{K})-d_2.
\end{eqnarray*}
Set
\[
A(f,m)=
\begin{pmatrix}
   \frac{\partial F_1}{\partial f} & \frac{\partial F_1}{\partial m} \\
   \frac{\partial F_2}{\partial f} & \frac{\partial F_2}{\partial m}
\end{pmatrix}
\]
Suppose $\lambda_1$ and $\lambda_2$ are two eigenvalues of $A(f,m)$. Then,
\bys
& & \lambda_1+\lambda_2=tr(A(f,m))=\frac{\beta}{2}\left[m(g_1-\frac{f}{K})+f(g_1-\frac{m}{K})\right]-d_1-d_2\\
& & \lambda_1\lambda_2=d_1d_2-\frac{\beta}{2}\left[f d_1(g_1-\frac{m}{K})+md_2(g_1-\frac{f}{K})\right]
\eys
A steady-steady solution $(f^*,m^*)$ of the system (4.1)-(4.2) is stable if and only if
\[ \lambda_1+\lambda_2<0, \h \lambda_1\lambda_2>0.\]

For Case 1, we have
\[
A(f,m)=A(0,0)=
\begin{pmatrix}
   -d_1 & 0 \\
   0 & -d_2
\end{pmatrix}
\]
This implies that the steady-state solution $(f,m)=(0,0)$ is linearly stable.
Moreover, ${\bar{ {\bf Z}}}=(0,0,0,0)$ is a global attractor.
\ \\
For Case 2, we have
\[
A(f_1,m_1)=A(0,0)=
\begin{pmatrix}
   -d_1 & 0 \\
   0 & -d_2
\end{pmatrix}
\]
\[
A(f_2,m_2)=A(\frac{4d_2}{\beta},\frac{4d_1}{\beta})=
\begin{pmatrix}
   -\frac{d_1d_2}{d_1+d_2}\ & \frac{d_2^2}{d_1+d_2}\\
   \frac{d_1^2}{d_1+d_2} & -\frac{d_1d_2}{d_1+d_2}
\end{pmatrix}
\]
Clearly, $(f_1,m_1)$ is asymptotically stable. For
$(f_2,m_2)$,  since
$det(A(f_2,m_2))=0$, it follows that
\[ \lambda_1=\frac{1}{2}tr(A)=-\frac{d_1d_2}{d_1+d_2}<0, \lambda_2=0.\]
This implies that $f(x,t)$ converges to $\bar{f}:=\frac{4d_2}{\beta}$ exponentially in $L^2(\Omega)$ as $t\rightarrow \infty$.
 However, $m(x,t)$ converges to either $0$ or $\frac{4d_1}{\beta}$ as $t\rightarrow \infty$. 
Therefore, the constant solution $(f_1,m_1)=(0,0)$ is asymptotically stable while $(f_2, m_2)=(\frac{4d_2}{\beta}, \frac{4d_1}{\beta})$ is unstable.

Now we discuss the stability for Case 3.

Clearly, $A(f_1,m_1)=A(0,0)$ is stable. For $(f_2, m_2)$,
\[ \lambda_1+\lambda_2<0 \Longleftrightarrow mg_1+fg_1< \frac{2mf}{K}+\frac{2(d_1+d_2)}{\beta},\]
which is true since all parameters are positive.

On the other hand,
\[ \lambda_1\lambda_2>0 \Longleftrightarrow \frac{\beta}{2}\left[fd_1(g_1-\frac{m}{K})+md_2(g_1-\frac{f}{K})\right]<d_1d_2,\]
which is equivalent to 
\[ \beta <\frac{8(d_1+d_2)}{K(1+b)^2}.\]
 This implies that the steady-state solution $(f_2,m_2)$ is asymptotically stable if
 \[ \beta <\frac{8(d_1+d_2)}{K(1+b)^2}\] and is unstable if 
\[ \beta >\frac{8(d_1+d_2)}{K(1+b)^2}\]
 However, for case 3, we require
 \[\beta>\frac{8(d_1+d_2)}{K}.\]
 Since $b>0$, no such $\beta$ satisfy both stability conditions. This implies that $(f_2,m_2)$ is always asymptotically unstable.
 
 Similarly, for $(f_3,m_3)$,
 \[ \lambda_1\lambda_2>0 \Longleftrightarrow \beta <\frac{8(d_1+d_2)}{K(1-b)^2}.\]
 This implies that the steady-state solution $(f_3,m_3)$ is asymptotically stable if
 \[ \beta <\frac{8(d_1+d_2)}{K(1-b)^2}\] and is unstable if 
\[ \beta >\frac{8(d_1+d_2)}{K(1-b)^2}\]
 
 Moreover, there exists a bifurcation of the solution $(f,m,s,r)$ with $t\rightarrow \infty$ when the birth rate $\beta$ changes across the critical value $\beta_0$.

\hfill Q.E.D.

\section{Conclusion}

In this paper, we studied a modified trojan Y Chromosome model. The modified model ensures that the solution to the reaction-diffusion system is nonnegative and uniformly bounded by the maximum carrying capacity $K$. The solution to the modified model makes better sense in the biological field. In addition, we analyzed the asymptotic behavior of the solution for the case where artificially introduced YY-male approaches $0$ as $t\rightarrow \infty$. It is shown that the natural female and male of the exotic species will be extinct if the birth rate is less than a critical number $\beta_0$, while they can survive for a long time if the birth rate is greater than the critical number $\beta_0$. Furthermore, we proved that there is a bifurcation phenomenon when the birth rate changes across the critical value $\beta_0=\frac{8(d_1+d_2)}{K}$.

\ \\
{\bf Acknowledgment.} 
The author would like to thank my colleague, Professor Xueying Wang, at Washington State University for showing me the interesting dynamics of the ODE model for trojan Y chromosome strategy.

\end{document}